\documentclass{aims}
\usepackage{amsmath}
  \usepackage{paralist}
 \usepackage[colorlinks=true]{hyperref}
\hypersetup{urlcolor=blue, citecolor=red}
\allowdisplaybreaks

\def\currentvolume{X}
 \def\currentissue{X}
  \def\currentyear{200X}
   \def\currentmonth{XX}
    \def\ppages{X--XX}

\newtheorem{theorem}{Theorem}[section]
\newtheorem{corollary}{Corollary}

\newtheorem{lemma}[theorem]{Lemma}

\theoremstyle{definition}

\title[Geophysics meets differential geometry] 
      {Geophysics and Stuart vortices on a sphere meet differential geometry}  

\author[{\L}ukasz Rudnicki]{}


 \email{lukasz.rudnicki@ug.edu.pl}

\thanks{{\L}ukasz Rudnicki is supported by the
Foundation for Polish Science (IRAP project, ICTQT, Contract No. 2018/MAB/5, cofinanced by the EU within the
Smart Growth Operational Programme)}

\begin{document}
\maketitle

\centerline{\scshape {\L}ukasz Rudnicki}
\medskip
{\footnotesize
 \centerline{International Centre for Theory of Quantum Technologies, University of Gda{\'n}sk}
   \centerline{80-308 Gda{\'n}sk, Poland}
} 

\bigskip

 \centerline{(Communicated by the associate editor name)}

\begin{abstract}
We prove new existence criteria relevant for the non-linear elliptic PDE of the form $\Delta_{S^2} u=C-he^{u}$, considered on a two dimensional sphere $S^2$, in the parameter regime $2\leq C<4$. We apply this result, as well as several previously known results valid when $C<2$, to discuss existence of solutions of a particular PDE modelling ocean surface currents.
\end{abstract}


Partial differential equations (PDEs) are one of the major modelling tools in physics, engineering and applied science in general. Looking from a perspective of a physicist, it is instructive to observe how distinct the type of interest in PDEs is, when it comes to the community working with them. For people subscribing themselves to the field of mathematics, it is usually of much more importance to discuss whether a given equation possesses a (well-behaved) solution or not. The question whether the solution can be found explicitly, or whether a numerical method serving such a purpose can be devised, is often of a lower importance. On the other hand, researchers describing physical systems usually tend to believe that since a given system exists in nature, the PDEs, which are supposed to model it, do have sufficiently regular solutions. Therefore, in this community it seems much more important to provide solutions and tools to obtain them, rather than to discuss the sole existence of the solutions.

The aim of this paper is to present an intermediate perspective, which uses a particular example relevant for geophysical applications (atmospheric flow on a planetary-scale, to be more precise) to bring deep results about existence of a particular non-linear elliptic PDE, related to differential geometry. Furthermore, these results return a feedback, which eventually helps fix the regimes of physically relevant parameters, in which the initial PDE can safely be considered.

We start a slightly retrodictive part of our discussion with the PDE
\begin{equation}
\Delta_{S^{2}}\psi-2\omega\cos\theta=\mathcal{F}\left(\psi\right),\label{Major}
\end{equation}
which, in certain circumstances, is supposed to govern ocean gyres \cite{Gyr1,Gyr2,Gyr3}. The stream function $\psi$, which depends on standard coordinates on a unit sphere
$\left(\theta,\varphi\right)$, where $0\leq\theta\leq \pi$ is the polar angle  and $0\leq\phi<2\pi$ is the angle
of longitude, encodes all information about the velocity field of an incompressible
flow. As usual

\begin{equation}\label{Lsph}
\Delta_{S^{2}}=\frac{\partial^{2}}{\partial\theta^{2}}+\cot\theta\frac{\partial}{\partial\theta}+\frac{1}{\sin^{2}\theta}\frac{\partial^{2}}{\partial\varphi^{2}},
\end{equation}
denotes the 2D Laplacian on the  sphere $S^2$. The parameter $\omega$
accounts for the effects due to rotation of the planet, while
$\mathcal{F}\left(\psi\right)$ is a free-to chose ``vorticity'' functional.

The form of $\mathcal{F}\left(\cdot\right)$ is not specified by the general
model under discussion, so that different settings (e.g. a linear function \cite{Gyr1}) are possible. A non-linear proposal in the form of 
\begin{equation}
\mathcal{F}\left(\psi\right)=ce^{d\psi}+g,\label{eq:F}
\end{equation}
with some real constants $c$, $d$ and $g$ seems to be of a particular relevance
\cite{Gyr2,Gyr3}, as it allows to enrich the model by the notion of \emph{Stuart
vortices}. Note that in \cite{Gyr2,Gyr3} the constants from (\ref{eq:F}) are
labelled with different letters (simply with $a,b,c$ respectively).
The notation used here will however prove itself to be suitable for
comparison with previous research. Consequently 
\begin{equation}
\Delta_{S^{2}}\psi-2\omega\cos\theta=ce^{d\psi}+g,\label{MajorM}
\end{equation}
is to be considered as our final model of ocean surface currents.

With the choice (\ref{eq:F}) we go back to 1967, when Stuart \cite{Stuart}
considered a variant of Eq. (\ref{Major}) on a plane (with $\Delta_{S^{2}}$
replaced by a standard, flat Laplacian in 2D), obviously with $\mathcal{F}\left(\cdot\right)$
given by (\ref{eq:F}), and with $\omega=0$. In original discussion by Stuart, the special case $g=0$ has been considered. Exact solutions of such
a steady two-dimensional Euler equation are nowadays commonly termed
as the Stuart vortices.

Much later, in 2004 Crowdy \cite{Crowdy} considered the same problem as Stuart had
done, however, with the underlying manifold being that of the unit sphere.
What follows from this single replacement is that the PDE considered
by Crowdy was of the form
\begin{equation}
\Delta_{S^{2}}\psi=ce^{d\psi}+g.\label{Major-1}
\end{equation}
As we can read in \cite{Crowdy} just after their ``(3.6)'' [which is our Eq. (\ref{Major-1})]:
\begin{center}
\emph{It will be shown in what follows that in the special case where
the parameters $d$ and $g$ are related by 
\[
g=\frac{2}{d},
\]
 then the general solution of (3.6) can be written in closed form.}
\par\end{center}

\noindent In other words
\begin{equation}
\Delta_{S^{2}}\psi=ce^{d\psi}+\frac{2}{d},\label{Major-2}
\end{equation}
which is a particular instance of (\ref{Major-1}), has been distinguished. In the discussion part the author further comments on that
fact, writing (see p. 398 in \cite{Crowdy}):
\begin{center}
\emph{It is intriguing that the (apparently special) choice $g=2/d$
has led both to the possibility of finding an explicit representation
for the general solution (...). While this might be coincidence, it
is easier to believe that the condition $g=2/d$ is a \textquoteleft solvability
condition\textquoteright{} for finding solutions (...)}
\par\end{center}

\noindent In the above historical detour we aimed at setting the scene
for the main considerations of this paper. While we are not going
to show that the choice $g=2/d$ is the conjectured solvability condition for
the problem (\ref{Major-1}) in a strict sense --- this equation admits a constant solution $\psi_{\textrm{const}}=d^{-1}\ln\left(-g/c\right)$ whenever $g/c<0$, which however is not of fluid-dynamical interest (vanishing flow) ---
in Sec. \ref{Sec2} we unravel the conundrum about why this choice is special. The gist lays with differential geometry. Furthermore, we scrutinize known statements, which pertain to solvability of a given family of elliptic PDEs, including the famous Kazdan and Warner "obstruction" criterion \cite{KzW}. In Sec. \ref{Sec3} we go back to geophysics, applying the tools summarized in Sec. \ref{Sec2} to the problem of solvability of (\ref{MajorM}). Then, in the most important Sec. \ref{Sec4}, we build on the results of Kazdan and Warner \cite{KzW} and Aubin \cite{AubinFr}, providing new \textit{sufficient conditions} for solvability of in a way a more general PDE
\begin{equation}
\Delta_{S^2} u=C-he^{u},\label{Elliptic}
\end{equation}
with $C=\mathrm{const}$ and $h(\theta,\phi)$ being a smooth function on the sphere, in the regime $2\leq C<4$. Note that all results quoted in Sec. \ref{Sec2} pertain to $C<2$. Finally, in Sec. \ref{Sec5} we apply this sufficient criterion to (\ref{MajorM}). We shall now summarize main results. 

Let: $F_{1} = \cos\theta$, $F_{2}=\sin\theta\cos\varphi$ and $F_{3}= \sin\theta\sin\varphi$ be the three standard, mutually orthogonal spherical harmonics of degree $1$ on $S^2$. Furthermore, let:
\begin{equation}
f_{i}=\boldsymbol{\nabla}h\cdot\boldsymbol{\nabla}F_{i}+\left(C-2\right)hF_{i},\qquad\qquad i\in\{1,2,3\},\label{Functions_fi}
\end{equation}
and
\begin{equation}
W_{ij}=\boldsymbol{\nabla}F_{i}\cdot\boldsymbol{\nabla}f_{j}+\left(C-2\right)F_{i}f_{j},\qquad\qquad i,j\in\{1,2,3\}.\label{MatrixW}
\end{equation}
Note that $f_i$ and $W_{ij}$ do depend on $C$. In Sec. \ref{Sec4} we prove the following theorem.
\begin{theorem}\label{Thmain1}
Let $h$ be positive somewhere. If for a given $2\leq C<4$:
\begin{itemize}
\item every $f_{i}$ defined in (\ref{Functions_fi}) changes the sign,
i.e. on $S^{2}$ assumes both positive and negative values;
\item $\sum_{i=1}^{3}\left|f_{i}\right|\geq\alpha$, with a gap $\alpha>0$;
\item $\det\left(\int_{S^{2}}d\Omega\Psi W\right)\neq0$ for all positive
functions $\Psi$, and the matrix $W$ defined in (\ref{MatrixW});
\end{itemize}
then there exists a smooth solution of (\ref{Elliptic}).
\end{theorem}
We supplement this theorem by the following statement, also proved in Sec. \ref{Sec4}.
\begin{theorem}\label{Thmain2}
If $W+W^{T}$ is either positive definite or negative definite for
every point on $S^{2}$, then $\det\left(\int_{S^{2}}d\Omega\Psi W\right)\neq0$
for all positive functions $\Psi$.
\end{theorem}
The existence results obtained for Eq. (\ref{MajorM}), which is here the major equation of interest in relation to its applications in geophysics, can be summarized as the following theorem which is a compilation of Corollaries \ref{Corollary1}-\ref{Corollary4} and Corollary \ref{Corollary7}.
\begin{theorem}
Let $c$, $d$, $g$ and $\omega$ be real constants in (\ref{MajorM}). Let $C= g d$ and $\varpi=\omega d$. If $\varpi\neq 0$ and:
\begin{itemize}
\item $C<0$, then there exists a solution of (\ref{MajorM}) if and only if $cd>0$;
\item $C=0$, then (\ref{MajorM}) does not have a solution;
\item $0<C<2$, then there exists a solution of (\ref{MajorM})  if and only if $cd<0$;
\item $C=2$, then (\ref{MajorM}) does not have a solution;
\item $2<C<13/6$, then there exists a solution of (\ref{MajorM})  if $cd<0$ and $\left|\varpi\right|<2\left(C-2\right)^{3/2}\left(9-4C\right)\sqrt{11-5C}$;
\item $13/6\leq C<4$, then there exists a solution of (\ref{MajorM})  if $cd<0$ and $\left|\varpi\right|<C-2$.
\end{itemize}
\end{theorem}

\section{Uniformization Theorem and Stuart vortices}\label{Sec2}
As the title of this chapter suggests, we are to recall a very well known
result from differential geometry. Given a two dimensional manifold
$\mathcal{M}$, we consider a metric tensor $\mathsf{g}_{0}$ on $\mathcal{M}$. Moreover,
we consider yet another metric tensor $\mathsf{g}$ on $\mathcal{M}$, which
can be expressed in terms of the former metric tensor as follows 
\begin{equation}
\mathsf{g}=e^{2\tilde u}\mathsf{g}_{0}.\label{conformal}
\end{equation}
The function $\tilde u$ is therefore just a conformal factor. More precisely,
we say that $\mathsf{g}$ and $\mathsf{g}_{0}$ are pointwise conformal.

Let $K_{0}$ be the Gaussian curvature of $\mathsf{g}_{0}$ and let
$K$ be a \emph{candidate} Gaussian curvature of $\mathsf{g}$. Furthermore,
let $\Delta_{0}$ be the Laplacian with respect to $\mathsf{g}_{0}$.
It is well known, and can easily be checked by a direct calculation,
that the function $\tilde u$ satisfies:
\begin{equation}
\Delta_{0}\tilde u=K_{0}-Ke^{2\tilde u}.\label{Diffgeo}
\end{equation}
However, on purpose we made $K$ to merely be a candidate for being
the Gaussian curvature, as \emph{a priori}, there is no guarantee
that (\ref{Diffgeo}) is solvable. For certain functions $K$ the
solution for the conformal factor might not exist \cite{KzW}.

The celebrated uniformization theorem states that if the Gaussian
curvature of $\mathsf{g}$ is constant, i.e. $K=\textrm{const}$, then
the metric $\mathsf{g}$ defined via (\ref{conformal}) exists. In
other words, there exists a function $\tilde u$ which satisfies (\ref{Diffgeo})
with $K=\textrm{const}$. 

\subsection{Connection with Stuart vortices on the sphere}
We are now ready to present the first observation reported in this
paper. If $\mathcal{M}=S^2$ and $\mathsf{g}_{0}$ is the metric tensor of a unit sphere,
we know that

\begin{equation}\nonumber
K_{0}=1\qquad\textrm{and }\qquad\Delta_{0}\equiv\Delta_{S^{2}},
\end{equation}
where the Laplacian $\Delta_{S^{2}}$ on the two dimensional sphere has been defined in (\ref{Lsph}). In such a case (\ref{Diffgeo}) becomes 
\begin{equation}\label{interm}
\Delta_{S^{2}}\tilde u=1-Ke^{2\tilde u}.
\end{equation}

If we then substitute 
\begin{equation}\nonumber
\tilde u=\frac{d}{2}\psi,
\end{equation}
which is nothing more than a rescaling of the function of interest,
and then multiply (\ref{interm}) by $2/d$, we get
\begin{equation}\nonumber
\Delta_{S^{2}}\psi=\left(-\frac{2K}{d}\right)e^{d\psi}+\frac{2}{d}.
\end{equation}
This equation is the same as (\ref{Major-2}), provided that one identifies $c=-2K/d$.

We can see that the very particular instance of the PDE (\ref{Major-1}), successfully solved by Crowdy \cite{Crowdy}, plays a very profound role in the core of the two-dimensional differential geometry. This fact partially explains why (\ref{Major-2}) is so special, so to speak,  hinting why its solutions were manageable. In fact, due to very recent results, we know even more. The major conclusion of \cite{ConstSol} can be summarized as follows (with slightly adjusted notation). Let $\lambda^*\geq 0$ be a minimum of a certain functional (see \cite{ConstSol} for the definition), an explicit form of which is of no relevance for the current discussion. Then one can prove that:
\begin{theorem}[Dolbeault, Esteban \& Jankowiak \cite{ConstSol}, p. 2: Theorem 1 and Corollary 2]\label{ThJan}
Let $\mathcal{M}$ be a smooth compact connected Riemannian manifold
of dimension $\dim\mathcal{M}=2$, without a boundary. If
$\lambda^{*}>0$ and $\tilde{u}$ is a smooth solution to
\[
-\frac{1}{2}\Delta_{0}\tilde{u}+\lambda=e^{\tilde{u}},
\]
then $\tilde{u}$ is a constant function if $0<\lambda<\lambda^{*}$.
If $\mathcal{M}=S^{2}$ and $\Delta_{0}=\Delta_{S^{2}}$, then the
choice $\lambda^{*}=1$ is optimal.
\end{theorem}

Clearly, under a substitution $\tilde{u}=d\psi+\ln\left(-cd/2\right)$
with $cd<0$, and after multiplying the above equation by $-2/d$,
we are back in Eq. (\ref{Major-1}), with $g=2\lambda/d$. From Theorem \ref{ThJan} we therefore learn that whenever $0<g<2/d$, the only smooth solution of (\ref{Major-1}) is a constant function. However, the authors of \cite{ConstSol} notice that
for $\lambda=\lambda^{*}$, i.e. $g=2/d$ in our case, the PDE in
question ``\emph{has non-constant solutions because of the
conformal invariance}''. We can see that finding non-trivial solutions
for $0<g<2/d$ had been impossible, therefore, only the very special
choice $g=2/d$ became fruitful.

Perhaps even more surprising result is to come. Since (\ref{Major-2}) is a particular instance of (\ref{Major-1}) --- we just need to set $g=2/d$ --- one can expect that a similar counterpart relevant for (\ref{MajorM}) will also be handy. However, in Sec. \ref{Sec3}  we will show that
\begin{corollary}\label{Corollary1}
If $\omega\neq 0$ then 
\begin{equation}\label{MajorMd}
\Delta_{S^{2}}\psi-2\omega\cos\theta=ce^{d\psi}+\frac{2}{d},
\end{equation}
does not have smooth solutions.
\end{corollary}
In other words, a mere and supposedly innocent inclusion of the rotation term, due to which  (\ref{Major-2}) is replaced by (\ref{MajorMd}) --- note that  (\ref{MajorMd}) reduces to (\ref{Major-2}) when $\omega=0$ --- has ultimately profound consequences. Before we establish this second result of our considerations, we shall first bring more background linking differential geometry with solvability of (\ref{Diffgeo}).

\subsection{Can all curvature functions $K$ be realized?}

In relation to (\ref{Diffgeo}), Kazdan and Warner \cite{KzW} in the very detail studied
the question whether a given function $K$ can be realized as the
Gaussian curvature of $\mathsf{g}$. As a matter of fact, their general
analysis also covers dimensions bigger than $2$ \cite{KzW,KzW1,KzW2}. To approach the
problem in a maximally broad sense, the authors abandoned geometric
considerations and focused on the following PDE
\begin{equation}
\Delta_{0}u=C-he^{u},\label{Diffgeo-1}
\end{equation}
with $C$ being a constant and $h$ being a smooth function. We shall exclude
the trivial case $h\equiv0$.
In this general setting, $\Delta_{0}$ remains to be the Laplacian
with respect to $\mathsf{g}_{0}$, and does not need to be equal to (\ref{Lsph}). 

In \cite{KzW}, the authors collected all the available information about
the solvability of (\ref{Diffgeo-1}), providing several supplementary
results which made some special cases more strict. Below, in a
series of lemmas we collect some of these findings --- those which
we find both the most informative and most useful for our discussion. 

Let $d\Omega$ denote the volume element, and let 
\begin{equation}\nonumber
\overline{h}=\frac{1}{\Omega\left(\mathcal{M}\right)}\int_{\mathcal{M}}d\Omega h,
\end{equation}
be the average value of $h$ on $\mathcal{M}$. By $\Omega\left(\mathcal{M}\right)$ we
denote the volume of $\mathcal{M}$ (area if $\textrm{dim}\mathcal{M}=2$). We have the
following \cite{KzW}:
\begin{lemma}\label{Lemma1}
If $C<0$ then $\overline{h}<0$ is a necessary (however not sufficient)
condition for solution of (\ref{Diffgeo-1}) to exist.
\end{lemma}
\begin{lemma}\label{Lemma2}
Solutions of (\ref{Diffgeo-1}) exist for all $C<0$ if and only if
$h\leq 0$ and $h$ is negative somewhere, so that the necessary condition from Lemma \ref{Lemma1} holds.
\end{lemma}
\begin{lemma}\label{Lemma3}
If $C=0$ and $\textrm{dim}\mathcal{M}=2$, then the solution of (\ref{Diffgeo-1})
exists if and only if $\overline{h}<0$ and $h$ is positive somewhere.
\end{lemma}
\begin{lemma}\label{Lemma4}
If $C>0$, $\textrm{dim}\mathcal{M}=2$ and $h$ is positive somewhere, there
is a constant $C_{+}\left(h,\mathcal{M}\right)>0$ such that the solution of
(\ref{Diffgeo-1}) exists if $C<C_{+}\left(h,\mathcal{M}\right)$. 
\end{lemma}
\begin{lemma}\label{Lemma5}
For the special case of the sphere, we have \cite{Moser}  $C_{+}\left(h,S^{2}\right)\geq2$.
\end{lemma}
Note that while Lemmas \ref{Lemma2} and \ref{Lemma3} are very precise (necessary and
sufficient conditions provided), the regime $C>0$ is quite vague.
While we know that the solutions exist for $0<C<C_{+}\left(h,\mathcal{M}\right)$,
we do not know what happens when $C\geq C_{+}\left(h,\mathcal{M}\right)$. 

\subsection{Integrable point singularities}
We shall stress that the above existence considerations hold under an assumption that no integrable point singularities on $\mathcal{M}$ are allowed. Technically speaking, if we integrate the left hand side of (\ref{Diffgeo-1}) over the entire manifold we get
\begin{equation}\label{int}
\int_{\mathcal{M}}d\Omega\Delta_{0}u=0,
\end{equation}
and consequently
\begin{equation}\label{GenConst}
\int_{\mathcal{M}}d\Omega he^{u}=\Omega\left(\mathcal{M}\right)C.
\end{equation}
This last condition is usually a supplementary tool linking the sign
of $C$ with the behavior of $h$, also playing a role of a constraint in variational derivation of  (\ref{Diffgeo-1}).

For example, the integrable point singularities on $S^2$ can be introduced by the terms
\begin{equation}\nonumber
u_{\textrm{sing}}^{\left(+\right)}\left(\theta,\varphi\right)=-C_{0}\ln\left[1+\cos\theta\right]\equiv C_{0}\ln\left[1+\tan^{2}\left(\frac{\theta}{2}\right)\right]-\ln2,
\end{equation}
or
\begin{equation}\nonumber
u_{\textrm{sing}}^{\left(-\right)}\left(\theta,\varphi\right)=-C_{0}\ln\left[1-\cos\theta\right]\equiv C_{0}\ln\left[1+\cot^{2}\left(\frac{\theta}{2}\right)\right]-\ln2.
\end{equation}
Note that in \cite{Crowdy} the second version --- this with the cot
function --- has been used, while for the sake of a more intuitive
interpretation we prefer here the variants with the cosinus. We can
clearly see that $u_{\textrm{sing}}^{\left(+\right)}$ has
a singularity at a south pole of the sphere, while $u_{\textrm{sing}}^{\left(-\right)}$
is singular at the north pole.

Calculating the Laplacian we can find that
\begin{equation}\nonumber
\Delta_{S^{2}}u_{\textrm{sing}}^{\left(\pm\right)}=C_{0},
\end{equation}
for both signs. Consequently, the integral similar to the one in (\ref{int})
does not vanish, giving the constant
\[
\int_{S^{2}}d\Omega\Delta_{S^{2}}u_{\textrm{sing}}^{\left(\pm\right)}=4\pi C_{0}.
\]
The volume element on the sphere as always reads $d\Omega=\sin\theta d\theta d\varphi$.

We can observe that if we allow such a point singularity and further decompose
\begin{equation}\nonumber
u=u_{\textrm{reg}}^{\left(\pm\right)}+u_{\textrm{sing}}^{\left(\pm\right)},
\end{equation}
where the regular part $u_{\textrm{reg}}^{\left(\pm\right)}$
obeys (\ref{int}), we find that
\begin{equation}\nonumber
\Delta_{S^{2}}u_{\textrm{reg}}^{\left(\pm\right)}=C-C_{0}-h_{\pm}e^{u_{\textrm{reg}}^{\left(\pm\right)}}.
\end{equation}
The modified candidate Gaussian curvature reads
\begin{equation}\nonumber
h_{\pm}=\frac{h}{\left(1\pm\cos\theta\right)^{C_{0}}}.
\end{equation}
We observe that the integrable point singularities allow one to shift
(increase or decrease) the constant $C$. The price to be paid is
that the curvature is either singular if $C_{0}>0$, or vanishes
on one of the poles when $C_{0}<0$.

Since we are potentially interested rather in decreasing the value of $C$ in the case when it exceeds $2$, we can see that the above method is not suitable for our purpose, as it will always lead to singular curvature functions. We will therefore not consider integrable point singularities from now on.

\subsection{Kazdan and Warner necessary condition}
In \cite{KzW} the following theorem has been proven.
\begin{theorem}[Kazdan \& Warner \cite{KzW}, p. 33: Theorem 8.8]
If $u$ is a solution to (\ref{Elliptic}), then \begin{equation}\label{KW}
    \int_{S^{2}}d\Omega e^{ u}\boldsymbol{\nabla}h\cdot\boldsymbol{\nabla}F=\left(2-C\right)\int_{S^{2}}d\Omega he^{u}F,
\end{equation}
for every $F$ being a spherical harmonics of degree 1, i.e. a solution of the equation  $\Delta_{S^{2}}F=-2F$. 
\end{theorem}
This result has been later on extended to $S^n$ \cite{KzW1}, and in many subtle ways improved. For a detailed account one can read a discussion in \cite{AubinBook} which starts on p. 233 with the proofs of (\ref{KW}), separately for $S^2$ and $S^n$ with $n>2$.

As already mentioned in the introduction, let:
\begin{equation}\label{SphHarm}
    F_{1} = \cos\theta, \qquad F_{2} = \sin\theta\cos\varphi, \qquad F_{3} = \sin\theta\sin\varphi,
\end{equation}
be the three standard, mutually orthogonal spherical harmonics of
degree $1$, which form a basis of the subspace labelled by the first
non-trivial eigenvalue of the Laplacian on the sphere. Clearly, every
$F$ which satisfies $\Delta_{S^{2}}F=-2F$ is of the form $F=\sum_{i=1}^{3}v_{i}F_{i}$,
with any numbers $\left\{ v_{i}\right\}$. Therefore, the Kazdan and
Warner criterion subsumes three independent conditions, each evaluated
for one of the $F_{i}$.

In other words, the functions (\ref{Functions_fi}) are the kernels of the Kazdan and Warner conditions, so that Eq.
(\ref{KW}) is equivalent to 
\begin{equation}\label{conconff}
\int_{S^{2}}d\Omega e^{u}f_{i}=0,\qquad\qquad i\in\{1,2,3\}.
\end{equation}
Let us also observe that the $3\times3$ matrix $W$ depending on the point
on $S^{2}$, defined in (\ref{MatrixW}), is nothing else than the $i$th Kazdan and Warner condition applied to $h=f_{j}$.

\section{Do solutions of  Eq. (\ref{MajorM}) even exist?}\label{Sec3}
In this section we are going to use the results presented above to study existence of solutions of our main equation of interest. To this end we use lemmas quoted in Sec. \ref{Sec2} to establish a chain of corollaries. As an initial step we use the substitution 
\begin{equation}\nonumber
\psi=\frac{1}{d}u-\omega\cos\theta,
\end{equation}
and multiply (\ref{MajorM}) by $d$, in order to transform it to the desired form of Eq. (\ref{Elliptic}) with $C=gd$, and the function $h$ explicitly given by
\begin{equation}\label{cadcur}
h_\omega=-cde^{-\omega d\cos\theta}.
\end{equation}
Interestingly, the candidate curvature (\ref{cadcur}), even though axially symmetric, is not antipodally symmetric. Its value on the north pole is $e^{-2\omega d}$ times smaller (or bigger) than the value it assumes on the south pole, depending on whether $\omega d$ is positive or negative. 

While considering the case $\omega\neq0$ we get the following:
\begin{corollary}\label{Corollary2}
If $gd<0$, Eq. (\ref{MajorM}) possesses a solution if and only if $cd>0$.
\end{corollary}
This is an immediate consequence of Lemma \ref{Lemma2} and the fact that (\ref{cadcur}) has a fixed sign over the entire sphere.
\begin{corollary} \label{Corollary3}
If $gd=0$, Eq. (\ref{MajorM}) does not have a solution.
\end{corollary}
We can see that $h_\omega$ defined in (\ref{cadcur}), if is positive somewhere, then
it must be positive everywhere on the sphere. Consequently, its average
is also positive, thus the assumptions on Lemma \ref{Lemma3} are not met.
\begin{corollary} \label{Corollary4}
If $0<gd<2$, Eq. (\ref{MajorM}) possesses a solution if and only if $cd<0$.
\end{corollary}
This time it is an immediate consequence of Lemmas  \ref{Lemma4} and  \ref{Lemma5}.
We note in passing that none of the above conditions depend on $\omega$, as long as $\omega\neq0$.

We remember that the choice $gd=2$ was of special importance in the case $\omega=0$, while we announced that for $\omega\neq0$ the corresponding problem does not have a solution. Now we are going to show Corollary \ref{Corollary1}, using the Kazdan and Warner necessary condition.

To this end, we shall evaluate the left hand side of (\ref{KW}), for $F=F_{1}$ and $h_\omega$
defined in (\ref{cadcur}). We get
\begin{equation}\label{non0}
\int_{S^{2}}d\Omega e^{u}\boldsymbol{\nabla}h_\omega\cdot\boldsymbol{\nabla}F_1=cd^{2}\omega\int_{S^{2}}d\Omega e^{u-\omega d\cos\theta}\sin^{2}\theta\neq0.
\end{equation}
Since for $gd=2$ also $C=2$, the right hand side of (\ref{KW}) vanishes.
Consequently, we find an example of Kazdan
and Warner obstruction. 

\section{New existence criteria for Eq. (\ref{Elliptic}) and $2\leq C<4$}\label{Sec4}
As we could see, existence of solutions of (\ref{Elliptic}) in the regime $C<2$
is well-understood. Let us stress that in line with the form of (\ref{Elliptic}) we currently restrict our main discussion
to $S^{2}$. On the other hand, quite little is known about solvability of (\ref{Elliptic}) when $C>2$.
We are going to partially fill this gap, proving sufficient existence
conditions for $2\leq C<4$, already presented in Theorem \ref{Thmain1} and Theorem \ref{Thmain2}. 

Before we prove both theorems, we need to quote and discuss a few
essential ingredients. We first translate from French (and adapt notation)
the following theorems by Aubin \cite{AubinFr}. Unfortunately these theorems
are not repeated in his very comprehensive book \cite{AubinBook} written in English.

\begin{theorem}[Aubin \cite{AubinFr}, p. 155: Théorème 4]\label{TAubin4}
On a compact Riemannian manifold $\mathcal{M}$ of dimension $\dim \mathcal{M}=n$, let 
\begin{equation}\nonumber
\mu_{n}=\left(n-1\right)^{n-1}n^{1-2n}\Omega\left(S^{n-1}\right)^{-1}.
\end{equation}
 Let $u\in H_{1}^{n}\left(\mathcal{M}\right)$ satisfy
\begin{equation}
\int_{\mathcal{M}}d\Omega u=0.\label{mean0}
\end{equation}
Then
\begin{equation}
\int_{\mathcal{M}}d\Omega e^{u}\leq\mathcal{C}\left(\mu\right)e^{\mu\left\Vert \boldsymbol{\nabla}u\right\Vert _{n}^{n}},\label{Onofri}
\end{equation}
where we can take $\mu=\mu_{n}+\epsilon$, with arbitrarily small
$\epsilon>0$, so that $\mu=\mu_{n}$ is the optimal choice for the
constant. 
\end{theorem}
\begin{theorem}[Aubin \cite{AubinFr}, p. 157: Théorème 6]\label{TAubin6}
On a compact Riemannian manifold $\mathcal{M}$ of dimension $\dim \mathcal{M}=n$, let $f_{i}$ for $i=1,\ldots,k$ be a family of differentiable functions,
such that every $f_{i}$ changes the sign and $\sum_{i=1}^{k}\left|f_{i}\right|\geq\alpha>0$.
Let $u\in H_{1}^{n}\left(\mathcal{M}\right)$ satisfy (\ref{mean0})
and
\begin{equation}
\int_{\mathcal{M}}d\Omega e^{u}f_{i}=0,\qquad\qquad i=1,\dots,k.\label{Constraints}
\end{equation}
Then (\ref{Onofri}) holds with $\mu=\mu_{n}/2+\epsilon$, with arbitrarily small
$\epsilon>0$, so that $\mu=\mu_{n}/2$ is the optimal choice for
the constant. 
\end{theorem}
The main conclusion of the second theorem is an improvement of the Moser-Trudinger-Onofri inequality (\ref{Onofri})
which occurs under additional constraints
(\ref{Constraints}). For more background on the Moser-Trudinger-Onofri inequality, see \cite{OnofriReview}. If $\dim\mathcal{M}=n=2$, then $\Omega\left(S^{1}\right)=2\pi$
is the circumference of a unit circle. Consequently, $\mu_{2}=1/16\pi$ is optimal in the generic case of $S^{2}$, while with the help of (\ref{Constraints}) it can be improved to $1/32\pi$.

Note that Theorem \ref{TAubin6} involves a positive constant $\alpha$, instead of just assuming $\sum_{i=1}^{k}\left|f_{i}\right|>0$. In Theorem \ref{Thmain1}, which inherits the assumptions from Theorem  \ref{TAubin6},  this constant is referred to as the gap, since it extensively plays such a role in the proof \cite{AubinFr} of Theorem \ref{TAubin6}.

As explained in many literature positions devoted to existence of
solutions to (\ref{Elliptic}), the inequality (\ref{Onofri}) can be
used to determine whether a minimum of the functional 
\begin{equation}\nonumber
J\left[u\right]=\int_{\mathcal{M}}d\Omega\left(\frac{1}{2}\left|\boldsymbol{\nabla}u\right|^{2}+Cu\right),
\end{equation}
is finite, while $u$ is subject to the constraint given by (\ref{GenConst})
for a suitable function $h$. Following \cite{KzW}, we write $u=v+\overline{u}$,
so that $\overline{v}=0$, and use (\ref{GenConst}) to find
\begin{equation}\nonumber
\overline{u}=\ln\left[\frac{C\Omega\left(\mathcal{M}\right)}{\int_{\mathcal{M}}d\Omega he^{v}}\right].
\end{equation}
If by $h_{\textrm{max }}<\infty$ we denote the supremum of $h$ on $\mathcal{M}$,
then
\begin{equation}\nonumber
\overline{u}\geq\ln\left[C\Omega\left(\mathcal{M}\right)/h_{\textrm{max }}\right]-\ln\left[\int_{\mathcal{M}}d\Omega e^{v}\right].
\end{equation}
Since $\boldsymbol{\nabla}u=\boldsymbol{\nabla}v$ and $v$ satisfies
(\ref{mean0}), we can apply Theorem \ref{TAubin4} (or Theorem \ref{TAubin6} if
additional conditions are satisfied) to get
\begin{equation}\nonumber
\overline{u}\geq\textrm{Const}-\mu\int_{\mathcal{M}}d\Omega\left|\boldsymbol{\nabla}v\right|^{2}.
\end{equation}
Consequently, we notice that the second term in $J\left[u\right]$
equals $C\Omega\left(\mathcal{M}\right)\overline{u}$, and because
$C>0$ we get
\begin{equation}\nonumber
J\left[u\right]\geq\left(\frac{1}{2}-C\Omega\left(\mathcal{M}\right)\mu\right)\int_{\mathcal{M}}d\Omega\left|\boldsymbol{\nabla}v\right|^{2}+\textrm{Const}.
\end{equation}
As a conclusion, if 
\begin{equation}
C<\frac{1}{2\Omega\left(\mathcal{M}\right)\mu},\label{BoundonC}
\end{equation}
the functional is bounded. For $S^{2}$, by virtue of Theorem \ref{TAubin4}
which gives $\mu=\mu_{2}$, we get $C<2$. From now on we shall again
restrict our attention to the two dimensional sphere.

In the realm of Theorem \ref{TAubin4} we just have a single constraint (\ref{GenConst}),
since (\ref{mean0}) does not directly affect $u$ in the discussion
presented above. Therefore, stationary points of $J\left[u\right]$
stem from varying ($\gamma$ is a Lagrange multiplier)
\begin{equation}\nonumber
J\left[u\right]-\gamma\int_{S^2}d\Omega he^{u},
\end{equation}
and satisfy
\begin{equation}\nonumber
\Delta_{S^{2}}u=C-\gamma he^{u}=0.
\end{equation}
The constraint (\ref{GenConst}) immediately yields $\gamma=1$, which means
that the minimum of $J\left[u\right]$ needs to be a solution of (\ref{Elliptic}).
Strict inequality in (\ref{BoundonC}) implies \cite{KzW, AubinBook} the existence
of this solution, however, we shall not elaborate here on that aspect.

It is now essential to check what are the implications of Theorem
\ref{TAubin6}. First of all, since $\mu=\mu_{2}/2$ we get the desired regime
$C<4$. Moreover, assuming that the requirements of Theorem \ref{TAubin6}
are met [for example, $f_{i}$ needs to change the sign so that
the set of functions satisfying the constraints (\ref{Constraints})
is not empty], the stationary points of $J\left[u\right]$ come
from
\begin{equation}\nonumber
J\left[u\right]-\gamma\int_{S^2}d\Omega he^{u}-\sum_{i=1}^{k}\beta_{i}\int_{S^2}d\Omega e^{u}f_{i},
\end{equation}
where $\left\{ \beta_{i}\right\} $ as well are Lagrange multipliers.
The function $u$ must satisfy (we already set $\gamma=1$)
\begin{equation}
\Delta_{S^{2}}u=C-\left(h-\sum_{i=1}^{k}\beta_{i}f_{i}\right)e^{u}=0.\label{NewEq}
\end{equation}
In other words, for $C<4$ and for some $\left\{ \beta_{i}\right\} $
to be specified with the help of (\ref{Constraints}), the above equation
has a solution. In fact, Aubin \cite{AubinFr} provided this result as Corollary 3 for $C=2$ and as a vaguely proved Corollary 4 for $C$ being
replaced by a function. In both cases Aubin set $k=3$ and took $f_{i}=F_{i}$, i.e. the basis of spherical harmonics of degree $1$ given in (\ref{SphHarm}). In these corollaries the primary candidate curvature $h$ is by virtue
of the variational principle replaced by  $h-F\left(h\right)$, where $F\left(h\right)$ is an
unspecified spherical harmonics of degree 1 (note that the notation in \cite{AubinFr} is different). A direct correspondence between these
results and (\ref{NewEq}) is immediate.

We are finally in position to prove Theorem \ref{Thmain1}.

\subsection{Proof of Theorem \ref{Thmain1}}

Let $2\leq C<4$. Since the functions $f_{i}$ defined
in (\ref{Functions_fi}) shall satisfy first two assumptions of the theorem for a given $C$, we know we can use these functions to apply Theorem \ref{TAubin6}. Therefore,
we know that in the discussed range of $C$, Eq. (\ref{NewEq}), with $k=3$ and $f_i$ defined in (\ref{Functions_fi}), has
a solution for some $\left\{ \beta_{1},\beta_{2},\beta_{3}\right\}$. To prove the theorem
under discussion, we need to show that necessarily $\beta_{1}=\beta_{2}=\beta_{3}=0$.

To this end, we apply the Kazdan and Warner criteria to the modified
curvature function
\[
h'=h-\sum_{i=1}^{3}\beta_{i}f_{i}\equiv h-\sum_{i=1}^{3}\beta_{i}\left[\boldsymbol{\nabla}h\cdot\boldsymbol{\nabla}F_{i}+\left(C-2\right)hF_{i}\right].
\]
As we already know that the solution of (\ref{Elliptic}) with $h$ replaced by $h'$
exists, we also know that this criteria must hold. Since the Kazdan
and Warner conditions for the ``bare'' function $h$, as shown in (\ref{conconff}), now play a
role of the constraints (\ref{Constraints}) selected for the purpose of Theorem \ref{TAubin6}, 
they vanish automatically. Therefore, we end up with three constraints
($i=1,2,3$)
\begin{equation}\nonumber
\sum_{j=1}^{3}\left(\int_{S^{2}}d\Omega e^{u}W_{ij}\right)\beta_{j}=0,
\end{equation}
with $W_{ij}$ defined in (\ref{MatrixW}). The solution of this equation is  trivial as desired, if and only if the matrix
$\int_{S^{2}}d\Omega e^{u}W_{ij}$ is non-singular. The last assumption
in Theorem \ref{Thmain1} assures that this is the case.

Finally, we shall observe that in the above method the solutions \textit{a priori} fulfill the constraints (\ref{conconff}). This fact does not diminish the validity or impact of the theorem, as all solutions (provided they exist) do anyway have to \textit{a posteriori} fulfill the Kazdan and Warner conditions.

\subsection{Proof of Theorem \ref{Thmain2}}
The conclusion of the second proposed theorem stems from the following
result in matrix theory.
\begin{theorem}[Horn \& Johnson \cite{Horn}, p. 510: Theorem 7.8.19 (Ostrowski-Taussky inequality); p. 511: Theorem 7.8.24]\label{ThforMatr}
Let $n\geq2$, let $H,K$ be two $n\times n$ Hermitian matrices and
let $A=H+iK$. If $H$ is positive definite then 
\begin{equation}\nonumber
\det H\leq\det H+\left|\det K\right|\leq\left|\det\left(H+iK\right)\right|=\left|\det A\right|.
\end{equation}
\end{theorem} 
Note that if $H$ is negative definite instead of being positive definite
then the same result can be applied for the matrix $-A$. Therefore,
we propose a slight generalization:
\begin{corollary}\label{CorroMat}
Let $n\geq2$, let $H,K$ be two $n\times n$ Hermitian matrices and
let $A=H+iK$. If $H$ is either positive definite or negative definite
then 
\begin{equation}\nonumber
0< \left|\det H\right|\leq\left|\det H\right|+\left|\det K\right|\leq\left|\det A\right|.
\end{equation}
\end{corollary}
An additional strict inequality on the left hand side, while is an obvious conclusion, will be important for our purpose. Interestingly, it is not easy and perhaps even not possible to abandon
the requirement of positive (negative) definiteness, in order to derive
a meaningful determinant inequality of a similar type. If one looks at the proof of
Theorem \ref{ThforMatr} (see \cite{Horn}), one observes that we need the identity $A\equiv H\left(I+iH^{-1}K\right)$ to hold, so
that $H$ needs to be invertible. However, a more important part is that the
matrix $H^{-1}K$ is diagonalizable and has only real eigenvalues
(see Corollary 7.6.2 in \cite{Horn}). This somehow less trivial fact stems
from Theorem 7.6.1 therein, which heavily uses the property of positive
definiteness of the matrix $H$.

While going back to the primary line of reasoning we need yet another standard
result from the theory of Hermitian matrices.
\begin{lemma}[For example in Horn \& Johnson \cite{Horn}, p. 239: as Theorem 4.3.1 (Weyl)]\label{Weyl}
Let $A,B$ be two Hermitian matrices of the same size and let $\lambda_{\min/\max}\left(A\right)$,
$\lambda_{\min/\max}\left(B\right)$ and $\lambda_{\min/\max}\left(A+B\right)$,
be the minimal/maximal eigenvalues of the matrices $A$, $B$ and
$A+B$ respectively. Then:
\begin{equation}\nonumber
\lambda_{\min}\left(A\right)+\lambda_{\min}\left(B\right)\leq\lambda_{\min}\left(A+B\right),
\end{equation}
\begin{equation}\nonumber
\lambda_{\max}\left(A+B\right)\leq\lambda_{\max}\left(A\right)+\lambda_{\max}\left(B\right).
\end{equation}
\end{lemma}
We obtain an immediate corollary.
\begin{corollary}
Let $A,B$ be two Hermitian matrices of the same size. If $A$ and
$B$ are both positive definite, $A+B$ is also positive definite.
If $A$ and $B$ are both negative definite, $A+B$ is also negative
definite.
\end{corollary}
We are ready to conclude the proof of Theorem \ref{Thmain2}. To this end,
we want to find sufficient conditions for $\det\left(\int_{S^{2}}d\Omega\Psi W\right)\neq0$
 to hold for all positive functions $\Psi$. Since $W$ is a real matrix (though
not necessarily symmetric), in order to apply the determinant inequalities
we need to set
\begin{equation}\label{HK}
H=\frac{1}{2}\int_{S^{2}}d\Omega\Psi\left(W+W^{T}\right),\qquad K=-\frac{i}{2}\int_{S^{2}}d\Omega\Psi\left(W-W^{T}\right).
\end{equation}
We can see from Corollary \ref{CorroMat} that if such $H$ is either positive definite or negative
definite, then $\left|\det\left(\int_{S^{2}}d\Omega\Psi W\right)\right|$
is bounded from below by a positive number. We thus just further need
to take care about positive/negative definiteness of $H$ defined in (\ref{HK}). Since, due to Lemma \ref{Weyl}, the
maximal eigenvalue is a convex function while the minimal eigenvalue
is concave (in fact they are both linear ``modulo'' min/max principle),
positive/negative definiteness of $H$ depends on whether $W+W^{T}$
has that property. Saying this we assume that the positive weight function $\Psi$ is well-behaved. For example
\begin{equation}\nonumber
\frac{1}{2}\int_{S^{2}}d\Omega\Psi\lambda_{\min}\left(W+W^{T}\right)\leq\lambda_{\min}\left(H\right).
\end{equation}
As in Theorem \ref{Thmain2} we in fact require that either $\lambda_{\min}\left(W+W^{T}\right)>0$ (positive definiteness)
or $\lambda_{\max}\left(W+W^{T}\right)<0$ (negative definiteness) for all points on the sphere $S^2$,
the proof is complete.

\section{New existence results for Eq. (\ref{MajorM})}\label{Sec5}

Since we have just established new tools allowing one to deal with
solvability of (\ref{Elliptic}) for $2\leq C<4$, we can now apply these techniques
to the primary problem discussed in this paper. Using the candidate curvature
$h_{\omega}$, given in (\ref{cadcur}), we can use (\ref{Functions_fi}) to compute:
\begin{align}
\begin{split}
f_{1}&=h_{\omega}\left[\left(C-2\right)\cos\theta-\omega d\sin^{2}\theta\right],\\
f_{2/3}& =h_{\omega}\left(C-2+\omega d\cos\theta\right)F_{2/3}.\label{FunctionsfOmega}
\end{split}
\end{align}
We immediately see that both $f_{2}$ and $f_{3}$ do change sign,
since $F_{2}$ and $F_{3}$ are $\varphi$-dependent. On the other
hand, for $C=2$ we find that $f_{1}$ is either positive or negative,
depending on the sign of $\omega c$ (note the constant $-cd$ in
front of $h_{\omega}$). We already know that in this case there are no solutions because the Kazdan and Warner criterion is violated [see (\ref{non0})]. Now, we can look at
the same from a different angle. For $C=2$ assumptions behind the
theorems proved in the previous section are not met. On the other
hand, for $C>2$, the function $f_{1}$ changes its sign, which can
easily be seen while comparing its values at north and south poles.

We can move on and verify the second property of this set of functions.
Since $\left|\cos\varphi\right|+\left|\cos\varphi\right|\geq1$, we
can see that
\begin{equation}\nonumber
\sum_{i=1}^{3}\left|f_{i}\right|\geq\left|h_{\omega}\right|\left|C-2\right|\left(\left|\cos\theta-\kappa\sin^{2}\theta\right|+2\left|\sin\theta\right|\left|1+\kappa\cos\theta\right|\right),
\end{equation}
where $\kappa=\omega d/\left(C-2\right)$. Since $h_{\omega}>0$,
we just need to assure ourselves that both terms in the sum on the
right hand side cannot simultaneously be $0$. The second term is
$0$ on both poles, where on the contrary the first term is equal
to $1$, or for $\theta=\arccos\left(-\kappa^{-1}\right)$, provided that $\left|\kappa\right|>1$.
However, in this special case, the first term is equal to $\left|\kappa\right|>0$. Being more precise, one can find the gap for this problem to be $\alpha=\left|c d\right|\left|C-2\right|e^{-\left|\omega d\right|}$. We leave the proof of that result as an exercise for interested readers.

We can see that the functions (\ref{FunctionsfOmega}) meet the criteria
posed by Theorem \ref{Thmain1} and inherited from Theorem \ref{TAubin6}. Therefore, sufficient conditions for solvability
of (\ref{MajorM}) will be determined by studying the corresponding matrix $W+W^{T}$.
While this matrix is rather cumbersome (so we do not provide it explicitly), we shall only concentrate here
on its eigenvalues. They are (we denote $\varpi=\omega d$)

\begin{equation}\nonumber
\lambda_{0}=2h_{\omega}\left(C-2+\varpi\cos\theta\right),\qquad\lambda_{\pm}=\frac{1}{4}h_{\omega}\left(\Theta\pm\sqrt{2\Delta}\right),
\end{equation}
where:
\begin{equation}\nonumber
\Theta=8+4C\left(C-3\right)+2\varpi^{2}+2\varpi\left(2\cos\theta-\varpi\cos2\theta\right),
\end{equation}
\begin{align}
\Delta = & 8\left(6-5C+C^{2}\right)^{2}+8\left[7+C\left(C-5\right)\right]\varpi^{2}+\varpi\left\{ 4\left[\varpi^{2}-4\left(C-3\right)\left(C-2\right)\right]\cos\theta\right.\nonumber \\
 &- \left.\varpi\left(4\left[\varpi^{2}+12+2C\left(C-5\right)\right]\cos2\theta+\varpi\left(4\cos3\theta-\varpi\cos4\theta\right)\right)\right\} +3\varpi^{4}.\nonumber
\end{align}
We shall look for the range of parameters $\varpi\in\mathbb{R}/\left\{ 0\right\} $
and $C\in\left(2,4\right)$ is which $W+W^{T}$ is either positive
or negative definite for all points on $S^{2}$ parameterized by two
angles $0\leq\varphi<2\pi$ and $0\leq\theta\leq\pi$. We observe
that even though the matrix in question does depend on $\varphi$,
its eigenvalues are axially symmetric. 

We begin our analysis with an observation that on the equator ($\theta=\pi/2$),
we get
\begin{equation}\nonumber
\Theta=8+4C\left(C-3\right)+4\varpi^{2}>8+4C\left(C-3\right)=4(C-1)(C-2)>0.
\end{equation}
This implies that there are always regions of the sphere where $\lambda_{+}$
is positive. This fact excludes an alternative of ``uniform'' negative
definiteness of $W+W^{T}$, so we just check whether all eigenvalues
are positive. Since $\lambda_{+}\geq\lambda_{-}$, we just need $\lambda_{-}>0$
and $\lambda_{0}>0$. The second requirement is quite easy, immediately
leading to a restriction
\begin{equation}
\left|\varpi\right|<C-2.\label{Firstrange}
\end{equation}

In the next step we shall show that, under the above condition, we
get $\Theta\geq0$ for all $0\leq\theta\leq\pi$. To this end we rewrite
$\Theta$ in terms of a new variable $X=\cos\theta$ which is in the
range $-1\leq X\leq1$
\begin{equation}\nonumber
\Theta\left(X\right)=8+4C\left(C-3\right)+4\varpi^{2}+4\varpi X-4\varpi^{2}X^{2}.
\end{equation}
Since this is a quadratic function with a negative coefficient in
front of $X^{2}$, such that it assumes a positive value for $X=0$,
to prove the claim it is enough to make sure that this function is
non-negative on the boundaries. We calculate
\begin{equation}\nonumber
\Theta\left(\pm1\right)=8+4C\left(C-3\right)\pm4\varpi.
\end{equation}
In the range (\ref{Firstrange}) we find
\begin{eqnarray}
\Theta\left(\pm1\right) & \geq & 8+4C\left(C-3\right)-4\left|\varpi\right|\nonumber \\
 & \geq & 8+4C\left(C-3\right)-4\left(C-2\right)\nonumber \\
 & = & 4\left(C-2\right)^{2}>0\nonumber.
\end{eqnarray}

Since we have shown that $\Theta\geq0$, the requirement $\lambda_{-}>0$
can be replaced by $\Theta^{2}-2\Delta>0$. We find that 
\begin{equation}
\Theta^{2}-2\Delta=64\left(C-2\right)^{2}\left[C-2+\varpi X\right]+16\varpi^{2}\left(4C-9\right)\left(1-X^{2}\right).\label{wyroznik}
\end{equation}
By virtue of (\ref{Firstrange}), whenever $4C-9\geq0$ this function
is always positive. We just need to check what happens if $C<9/4$.
Since $\Theta^{2}-2\Delta$ again is a quadratic function manifestly positive on the boundaries (when $|X|=1$),
we only need to test the minimum of this function in its domain. The
minimum occurs at 
\begin{equation}\nonumber
X_{\min}=-2\frac{\left(C-2\right)^{2}}{\varpi\left(9-4C\right)}.
\end{equation}
We can either allow $\left|X_{\min}\right|\geq1,$ or accept $\left|X_{\min}\right|<1$,
provided that (\ref{wyroznik}) is positive at $X_{\min}$. The first
option renders
\begin{equation}
2\frac{\left(C-2\right)^{2}}{\left(9-4C\right)}\geq\left|\varpi\right|,\label{range1}
\end{equation}
while the second one gives
\begin{equation}
2\frac{\left(C-2\right)^{2}}{\left(9-4C\right)}<\left|\varpi\right|<2\frac{\left(C-2\right)^{3/2}\sqrt{11-5C}}{\left(9-4C\right)},\label{range2}
\end{equation}
but only if $C\leq11/5$. In fact the range (\ref{range2}) is not empty
when $C\leq13/6<11/5$. Surprisingly, only till that value of $C$, namely $C=13/6$, the condition
(\ref{range1}) gives us more than (\ref{Firstrange}). In other words,
the final conclusion simplifies quite a bit. We can summarize it in
the following corollary.
\begin{corollary}\label{Corollary7}
Let $C=gd$ and $\varpi=\omega d$. For $\varpi\neq 0$ there exists a solution of (\ref{MajorM}) if $cd<0$ and
\begin{equation}\nonumber
\left|\varpi\right|<\begin{cases}
2\left(C-2\right)^{3/2}\left(9-4C\right)\sqrt{11-5C} & \quad\textrm{when }\; 2<C\leq13/6\\
C-2 & \quad\textrm{when }\; 13/6\leq C<4
\end{cases}.
\end{equation}
\end{corollary}

\section{Discussion}
The main result reported in this paper is Theorem \ref{Thmain1}. While one can think, for example reading the case study presented in Sec. \ref{Sec5}, that it might be difficult to get exact conclusions in scenarios more complex than Eq. (\ref{MajorM}), the theorem is also a handy tool for numerical studies. Since we work with $S^2$, all the conditions listed in Theorems \ref{Thmain1} and \ref{Thmain2} can approximately be checked numerically. As the matrix $W$ is of dimension $3$, its eigenvalues, no matter how cumbersome, can explicitly be found. Then for fixed values of the involved parameters, checking positive or negative definiteness of this matrix is an optimization problem on $S^2$. 

Diverging a bit from the main geophysics-oriented theme, we can realize that the results presented here can relatively easily be generalized. First of all, Kazdan and Warner criteria have been extended to cover $n$-dimensional spheres \cite{KzW1}. Also theorems by Aubin quoted in Sec. \ref{Sec4} do naturally cover such a case via the constant $\mu_n$. However, we leave the problem of generalizing Theorem \ref{Thmain1} to more dimensions for the future, simply because $S^n$ with $n>2$ is not relevant for geophysics applications. 

Moreover, the criterion by Kazdan and Warner appeared to be significantly simpler for the standard Liouville equation in $\mathbb{R}^2$ \cite{KzWR2}. We believe a variant of Theorem \ref{Thmain1} relevant for this scenario can also be of interest. The same remark applies to a similar problem relevant for the torus \cite{Sakajo}. Even though a counterpart of the parameter $g$, and therefore also $C$, is not constant in this scenario, the results by Aubin (see Corollaire 4 in \cite{AubinFr}) provide a good starting point for future research in this direction.

Last but not least, we notice that a variant of Theorem \ref{TAubin6} has just been proven \cite{Impggg}, where the constraints are specified to cover several subspaces of the $n$ dimensional Laplacian on the sphere. In this case the optimal constant $\mu$ is divided by a number depending on the degree of a particular spherical t-design on the sphere \cite{tdesigns}. However, the obstacle in using this result to extend the range of the parameter $C$ is that the Kazdan and Warner criteria do only hold in relation to the first non-trivial subspace of the Laplacian. Beyond spherical harmonics of degree 1, there is an additional term which depends on $|\boldsymbol{\nabla} u|$, and consequently spoils all the argumentation. On the other hand, a handful of geometric generalizations of the Kazdan and Warner criteria (see eg. \cite{exte}) have been developed over the years, therefore, several promising possibilities are still there.

\section*{Acknowledgments} I thank Darren Crowdy for fruitful discussions and Jerry Kazdan for pointing out Ref. \cite{exte}.









\medskip
Received xxxx 20xx; revised xxxx 20xx; early access xxxx 20xx.
\medskip

\end{document}